\documentclass[11pt]{amsart}

\usepackage{geometry}
\usepackage[all]{xy}

\usepackage{amssymb}
\usepackage{calrsfs}
\usepackage[mathscr]{eucal}
\usepackage[T1]{fontenc} 
\usepackage{textcomp}
\usepackage{times}
 \usepackage[scaled=0.92]{helvet} 
\renewcommand{\tilde}{\widetilde}
\renewcommand{\bar}{\overline}
\newcommand{\isomto}{\overset{\sim}{\rightarrow}}
\newcommand{\ab}{\mbox{{\tiny \textup{ab}}}}
\DeclareMathOperator{\F}{\mathbf{F}}
\DeclareMathOperator{\Z}{\mathbf{Z}}
\DeclareMathOperator{\A}{\mathbf{A}}
\DeclareMathOperator{\cA}{\mathcal{A}}
\DeclareMathOperator{\cO}{\mathcal{O}}
\DeclareMathOperator{\fp}{\mathfrak{p}}
\DeclareMathOperator{\q}{\mathfrak{q}}
\DeclareMathOperator{\ff}{\mathfrak{f}}
\DeclareMathOperator{\fn}{\mathfrak{n}}
\DeclareMathOperator{\fm}{\mathfrak{m}}
\usepackage{bbm}
\DeclareMathOperator{\one}{\mathbbm{1}}
\DeclareMathOperator{\ord}{\mathrm{ord}}
\DeclareMathOperator{\Places}{\mathrm{Places}}
\DeclareMathOperator{\Aut}{\mathrm{Aut}}
\DeclareMathOperator{\rec}{\mathrm{rec}}

\newtheorem{theorem}{Theorem}
\newtheorem{lemma}[theorem]{Lemma}
\newtheorem{corollary}[theorem]{Corollary}

\begin{document}

\date{\today\ (version 1.0)} 
\title{Curves, dynamical systems and weighted point counting}
\author[G.~Cornelissen]{Gunther Cornelissen}
\address{\normalfont Mathematisch Instituut, Universiteit Utrecht, Postbus 80.010, 3508 TA Utrecht, Nederland}
\email{g.cornelissen@uu.nl}


\begin{abstract} \noindent Suppose $X$ is a (smooth projective irreducible algebraic) curve over a finite field $k$. Counting the number of points on $X$ over all finite field extensions of $k$ will not determine the curve uniquely. Actually, a famous theorem of Tate implies that two such curves over $k$ have the same zeta function (i.e., the same number of points over all extensions of $k$) if and only if their corresponding Jacobians are isogenous. 
We remedy this situation by showing that if, instead of just the zeta function, \emph{all} Dirichlet $L$-series of the two curves are equal via an isomorphism of their Dirichlet character groups, then the curves are isomorphic up to ``Frobenius twists'', i.e., up to automorphisms of the ground field. 
Since $L$-series count points on a curve in a ``weighted'' way, we see that weighted point counting determines a curve. In a sense, the result solves the analogue of the isospectrality problem for curves over finite fields (also know as the ``arithmetic equivalence problem''): it says that a curve is determined by ``spectral'' data, namely, eigenvalues of the Frobenius operator of $k$ acting on the cohomology groups of all $\ell$-adic sheaves corresponding to Dirichlet characters. 
The method of proof is to shown that this is equivalent to the respective class field theories of the curves being isomorphic as dynamical systems, in a sense that we make precise. 
  \end{abstract}

\maketitle

\section{Introduction}

Let $X$ denote a smooth projective curve over a finite field $k=\F_q$ of characteristic $p$. Its zeta function encodes precisely the sequence of integers given by the number of points of $X$ over finite extensions $\F_{q^n}$ of $k$. As observed by Turner \cite{Turner}, a celebrated theorem of Tate (\cite{Tate}, Thm.\ 1) implies that an equality of the zeta functions of two curves over $k$ (so, equality of their number of points over all extensions of $k$) is equivalent to the curves having isogenous Jacobians over the field $k$. Since there are plenty of non-isomorphic curves with isogenous Jacobians, this shows that \emph{the $k$-isomorphism type of a curve over a finite field $k$ is not determined solely by point counting over field extensions.} For a very explicit example, consider the curves of genus two given by $$ X_{\pm} \colon y^2=x^5\pm x^3+x^2-x-1 \mbox{ over }k = \F_{3}.$$ It is know that these curves are not (even geometrically) isomorphic, but have the same zeta function, and furthermore their Jacobian is absolutely simple \cite{Howe}. 

The aim of the current paper is to show that if this counting is extended from using only the zeta function to abelian $L$-series, then this \emph{does} determine the curve uniquely, as follows:  
\begin{theorem} \label{cor}
Let $X,Y$ denote two smooth projective curves over a finite field $k$. Assume that there is an isomorphism $\psi \colon G_{k(X)}^{\ab} \isomto G_{k(Y)}^{\ab}$ of Galois groups of the maximal abelian extensions of their respective function fields,  such that the corresponding $L$-series match: 
$$ L(X,\chi) = L(Y,(\psi^{-1})^{\vee}\chi) $$
for all characters $\chi \in G_{k(X)}^{{\ab} \vee}$. Then $X$ is isomorphic to a \emph{Frobenius twist} $Y^\sigma$ of $Y$ over $k$, where $\sigma$ is an automorphism of the ground field $k$.  

For the conclusion to hold, it suffices that the $L$-series of geometric characters match (i.e., characters that factor through a finite geometric extension of the function field).  
\end{theorem}
We have the following ``explicit'' counting form of the $L$-series as a function of $T=q^{-s}$:  
\begin{align*}
L(X&,\chi,s) = \\ \exp & \left( \sum_{n \geq 1} \frac{T^n}{n} \sum_{d \mid n} \sum_{\zeta \in \mu_\infty} \zeta^{n/d}  \sum_{m \mid d} \mu\left(\frac{d}{m}\right) N(X,m,\zeta,\chi)    \right), 
\end{align*}
where $\mu$ is the M\"obius function, $\mu_\infty$ is the set of complex roots of unity, and $N(X,m,\zeta,\chi)$ is the number of points $x \in X(\F_{q^{m}})$ such that any corresponding place $[x]$ of $\F_q(X)$ is unramified for $\chi$, and $\chi(\mathrm{Frob}_{[x]}) = \zeta.$ 
Thus, some kind of \emph{weighted point counting} will determine the curve up to isomorphism, modulo so-called \emph{Frobenius twists} $Y^\sigma$ of $Y$, that are induced by an automorphism $\sigma \in \Aut (k)$ of the ground field $k$. Note that if the ground field $k$ is $\F_p$ for a prime $p$, then there are no such non-trivial twists. In the above example over $\F_{3}$, we now know that point counting weighted by characters does determine the isomorphism type of the curve. If $k=\F_{p^d}$ for $d>1$ however, then there exist curves that are not isomorphic to their twists, not even over the separable closure $\bar{k}$ (for example, an elliptic curve $E/k$ whose $j$-invariant $j(E)$ is not invariant under $\sigma \in \Aut (k)$ will have $E^\sigma$ not isomorphic to $E$ over $\bar{k}$). The presence of Frobenius twists is natural, since the objects we consider depend only on the function fields \emph{as fields} (not as extensions of the ground field), whereas morphisms between curves correspond to function field isomorphisms that \emph{fix} the ground field. 

The above theorem is a mere corollary of the main result. The application arises from associating to a curve over a finite field a certain noncommutative dynamical system, that derives from its abelian class field theory. More precisely, consider the topological space 
$$ S_X:=G_{{k}(X)}^{\ab} \times_{\cO^*_X} \cO_X, $$
(where $\cO_X$ is the ring of integral adeles of $X$),% and $\bar k$ is a separable closure of $k$), 
which is the quotient of the space $G_{{k}(X)}^{\ab} \times \cO_X$ by the equivalence $$(\gamma,\rho) \sim (\gamma \rec_X(u)^{-1}, u\rho) \mbox{ for any } u \in\cO_X^*,$$ where $\rec_X \colon \A^*_X \rightarrow G_{k(X)}^{\ab}$ is the Artin reciprocity map, with $\A^*_X$ the ideles of $X$. On this space, we have an action of the semigroup $I_X$ generated by the places of $k(X)$, given by $\fn \in I_X$ acting as 
$$ \fn \ast [(\gamma, \rho)] := [(\gamma \cdot \rec_X(s(\fn))^{-1}, s(\fn) \rho)], $$
where $s$ is any section of the natural map $\A^*_X  \rightarrow I_X$ (which from now on we leave out of the notation). 

The intuition behind the introduction of this space is the following: one would like to let ideals act by multiplication with their associate Frobenius elements in the Galois group, but those are only defined ``up to inertia'', and this is compensated for by the equivalence relation. 

An \emph{isomorphism} (a.k.a.\ topological conjugacy) of such systems is a pair $(\Phi,\varphi)$ consisting of a homeomorphism $\Phi \colon S_X \isomto S_Y$ and a semigroup isomorphism $\varphi \colon I_X \rightarrow I_Y$ such that $\Phi(\fn \ast x) = \varphi(\fn)\ast \Phi(x) \mbox{ for all }x \in S_X\mbox{ and }\fn \in I_X.$ We say $(\Phi,\varphi)$ is \emph{degree-preserving} if $\deg(\varphi(\fn)) = \deg(\fn)$ for all $\fn \in I_X$.

The main result is then
\begin{theorem} \label{main}
Let $X,Y$ denote two smooth projective curves over a finite field $k$. Assume that there is a degree-preserving isomorphism (i.e., topological conjugacy) of the associated dynamical systems $ (S_X,I_X) \isomto (S_Y,I_Y); $ then $X$ and $Y$ are isomorphic over $k$, up to Frobenius twists. 
\end{theorem}

In \cite{CM}, one finds the analogue of these theorems for number fields. The way to deduce Theorem \ref{cor} from Theorem \ref{main} is identical to that in \cite{CM}, and will be sketched in the next section; the part about restricting to geometric characters will be proven in the final section. As will be explained in the middle two sections of the paper, the proof of Theorem \ref{main} in the case of function fields will be somewhat different from that in number fields, and rely on Theorem \ref{thmuch} below---which might be of independent interest---that we deduce from the proof of Uchida's anabelian theorem for global function fields in \cite{Uchida}. 

The dynamical system $(S_X,I_X)$ can in its turn be encoded in a non-involutive Banach algebra $\cA_X$, and (as in the work of Davidson and Katsoulis \cite{DK}) isomorphism of these algebras will imply so-called \emph{piecewise conjugacy} of the associated dynamical systems. One may then show, as in \cite{CM}, Section 6, that degree preservation for this system (and ergodicity w.r.t.\ the natural measure induced from Haar measure, as established in \cite{LLN}) implies full topological conjugacy. Hence the main theorem could be rephrased in operator theoretical language as follows: a degree-preserving isomorphism of $\cA_X \isomto \cA_Y$ implies an isomorphism of the curves $X$ and $Y$ up to Frobenius twists. Degree-preservation may further be coded into preservation of a one-parameter subgroup of automorphisms of the algebra, leading to a description in terms of \emph{quantum statistical mechanical systems}. In this paper, we will not consider this ``physical'' viewpoint any further, but it was the original source for the discovery of Theorem \ref{main}. 

We do remark that our dynamical system is a function field cousin of the system that Ha and Paugam (\cite{HP}) introduced for number fields, as a generalization of the famous Bost-Connes dynamical system \cite{BC}, that was the first source to describe  abelian class field theory for the field of rational numbers by a dynamical system. Our system was recently also considered by Neshveyev and Rustad \cite{NR}. Other dynamical systems of a similar flavor have been associated to global function fields, in the work of Jacob (\cite{Jacob}) and Consani and Marcolli (\cite{ConsMar}), using Hayes' explicit class field theory for function fields from the theory of Drinfeld modules, but these systems seem to serve other purposes than the one of the current paper. 

The analogue of Theorem \ref{cor} for number fields has been proven without reference to dynamical systems and strengthened to an effective statement by Bart de Smit (unpublished), 
 using ideas from the theory of arithmetic equivalence. The effective statement is that every number field admits a character of order 3, whose $L$-series does not occur as $L$-series for any other character on any other number field. The method of proof of this result does not generalize readily to global function fields (compare with the discussion in \cite{CKZ} on arithmetic equivalence in function fields), so it is an interesting question to investigate whether or not Theorem \ref{cor} can similarly be made effective. It seems that the method of proof in characteristic zero (based on representation theory) breaks down in positive characteristic. Furthermore, as rational functions of $q^{-s}$, function field $L$-series seem to carry less information than their number field cousins (compare with the discussion in Gross \cite{Gross}).  
 
 Nevertheless, one may wonder whether or not a ``smaller'' invariant constructed from abelian class field theory determines a curve. In the example of the curves $X_{\pm}$ above, we used the explicit class field theory implemented in Magma \cite{Magma} to find the following. Both curves have exactly four places of degree two. Consider a divisor of the form $D=2P+Q+R$ with $P,Q,R$ three such different places of degree two, and let $S$ denote the remaining place. There are 12 different choices for such $D$. Amongst the abelian covers of degree 3 corresponding to the ray class group of $D$ and in which $S$ is totally split, we find the following distinction between the two curves: for every choice of $S$ on $X_-$, there is such a cover with no $k$-rational points, whereas for exactly two choices of $S$ on $X_+$, all such covers have $k$-rational points (so in particular, the list of zeta functions of these covers distinguishes the original curves). 
   
  A related interesting question is whether an equality between the set of zeta functions of all abelian covers of a curve (i.e., the isogeny classes of Jacobians of all abelian covers) determines an isomorphism of curves, up to Frobenius twist. This looks like a more geometrical question. 

The first theorem provides a new look at the problem of arithmetic equivalence, which is an analog of the isospectrality problem for Riemannian manifolds, cf.\ \cite{CdJ} for more on a Riemannian analog of Theorem \ref{cor}. The theorem can be reformulated as saying that the ``spectral invariants'' given by \emph{the characteristic polynomials of the $k$-Frobenius acting on the cohomology of the $\ell$-adic sheaf corresponding to the character $\chi$ determine the curve up to Frobenius twists.} 

The second theorem complements the program of anabelian geometry and/or Torelli theorems, trying to characterize geometrical objects by topological and/or group theoretical constructs. We see that a noncommutative topological dynamical system built from the \emph{abelian} class field theory reconstructs the field up to isomorphism (viz., the curve up to Frobenius twists). A recent result by Bogomolov, Korotiaev and Tschinkel \cite{BKT} and Zilber \cite{Zilber} that is also of an ``abelian anabelian'' nature, might be related: a curve of genus $\geq 2$ over a finite field is determined up to geometric isomorphism and Frobenius twists by the data of its Jacobian as a \emph{group}, the degree-one cycles as homogeneous space for this group, with its natural embedding of the curve. I do not know of an exact formal relation between their result and the results of this paper.  

\section{From Theorem \ref{main} to Theorem \ref{cor}} \label{buh}

\begin{theorem} \label{main3}
If $X$ and $Y$ satisfy conditions of Theorem \ref{cor}, then there is a degree-preserving isomorphism $\varphi \colon I_X \isomto I_Y$ such that for every finite abelian extension $K'=\left(k(X)^{\ab}\right)^N$ of $k(X)$ (with $N$ a subgroup in $G_{k(X)}^{\ab}$) and every place $\fp$ of $k(X)$ unramified in $K'$, the place $\varphi(\fp)$ is unramified in the corresponding field extension $L':=\left(k(Y)^{\ab}\right)^{\psi(N)}$ of $k(Y),$ and the induced map $\psi \colon \mathrm{Gal}(K'/k(X)) \rightarrow \mathrm{Gal}(L'/k(Y))$ satisfies
$\psi \left( \mathrm{Frob}_{\fp} \right) = \mathrm{Frob}_{\varphi(\fp)}. $ 
\end{theorem}

\noindent \emph{Proof.}\ 
Represent the $L$-series as a sum over effective divisors and use literally the same proof as for the number field case in \cite{CM}, Theorem 12.1. (The result there is formulated in terms of characters, but translates one-to-one to the statement here.) $\Box$

\begin{corollary}
If $X$ and $Y$ satisfy conditions of Theorem \ref{cor}, then there is a degree-preserving isomorphism of dynamical systems $(S_X,I_X) \isomto (S_Y,I_Y)$. 
\end{corollary}

\noindent \emph{Proof.}\ 
As in Corollary 5.5 of \cite{CM}, $S_X$ has a dense subspace given by classes $[(\gamma, \rho)]$ where none of the entries of $\rho $ is zero. For such an element, we define $$\Phi([(\gamma,\rho)]) = [(\psi(\gamma), \tilde{\rho} \fm_\rho)],$$ where we write $\rho = \rho' \cdot \fm_\rho$ for some unit idele $\rho' \in \cO_X^*$ and place $\fm_\rho \in I_X$, and $\tilde\rho \in \cO_Y^*$ is defined by $\psi(\rec_X(\rho')) = \rec_Y(\tilde \rho)$. We have used the previous theorem, that implies that $\psi$ maps the ``ramification group'' $\cO_X^*$ to the corresponding $\cO_Y^*$. 

This definition is compatible with the action of places as follows: for $\fn \in I_X$, \begin{align*} & \Phi(\fn \ast [(\gamma,\rho)]) = \Phi([(\gamma \rec_X(\fn)^{-1}, \fn \rho)]) \\ & = [(\psi(\gamma) \psi(\rec_X(\fn))^{-1}, \tilde{\rho} \varphi(\fm \fn))] \\ &=  [(\psi(\gamma) \rec_Y(\varphi(\fn))^{-1}, \tilde{\rho} \varphi(\fm)\varphi(\fn))] \\ & = \varphi(\fn) \ast [(\psi(\gamma), \tilde \rho \varphi(\fm))] = \varphi(\fn) \ast \Phi([\gamma,\rho)]).\end{align*} Here, we have used the previous theorem in replacing $\psi(\rec_X(\fn))$ by $\rec_Y(\varphi(\fn))$. It is easy to see that $\Phi$ extends to a homeomorphism on all of $S_X$ that remains compatible with the actions. 
$\Box$

We now continue with the proof of Theorem \ref{main}. 

\section{Revisiting the anabelian theorem of Uchida}

\begin{theorem} \label{thmuch}
Suppose that $K$ and $L$ are global function fields, such that the following exist:
\begin{enumerate}
\item[\textup{(a)}] a bijection of places $\Places(K) \overset{\varphi}{\rightarrow} \Places(L)$;
\item[\textup{(b)}] a group isomorphism $G_K^{\ab} \overset{\Phi}{\rightarrow} G_L^{\ab}$; 
 \item[\textup{(c)}] for any place $v \in \Places(K)$, an isomorphism of the multiplicative groups of the corresponding local fields $$ K_v^* \overset{\Phi_v}{\rightarrow} L_{\Phi(v)}^*$$ that is valuation preserving, i.e., such that $$\ord_{\varphi(v)}(\Phi_v(x)) = \ord_v(x) \mbox{ for all }x \in K_v;$$
\end{enumerate}
with the property that the following diagram commutes: 
\begin{equation*} \label{comm}  \xymatrix{ \A_K^* \ar@{->}[r]^{\rec_K} \ar@{->}[d]^{\prod \Phi_v} & G_K^{\ab}\ar@{->}[d]^{\Phi}   \\   \A_L^* \ar@{->}[r]^{\rec_L}  & G_L^{\ab}  }, \end{equation*} 
where $\A_K^*$ is the idele group of $K$ and $\rec_K$ is the Artin reciprocity map. 
Then the fields $K$ and $L$ are isomorphic (as fields), and the corresponding smooth curves are isomorphic up to Frobenius twists. 
\end{theorem} 

\noindent \emph{Proof.}\ The proof is a combination of some intermediate results in Uchida's proof of the anabelian theorem for global function fields \cite{Uchida}, and we now outline how to perform this combination, quoting the relevant results from \cite{Uchida}. 

Since $K^*$ is the kernel of the reciprocity map, from the assumptions we find a multiplicative map of fields
$$ K^* \isomto \ker(\rec_K) \overset{\prod \Phi_v}{\longrightarrow} \ker(\rec_L) \isomto L^*, $$
which we extend by $0_K \mapsto 0_L$. We will show that the resulting map is additive. For that, we first observe that condition (c) implies that $\Phi_v$ induces a multiplicative isomorphism of residue fields (since it preserves valuations). The claim is that this reduced map is also additive (analogue of Lemma 10 in \cite{Uchida}). Since a multiplicative map that is additive modulo all places is additive itself, this proves that $K$ and $L$ are isomorphic as fields. 

We now turn to proving additivity modulo all places. We first quote Lemma 8 from \cite{Uchida}, which shows that the given assumptions imply that the constructed map $K^* \rightarrow L^*$ has the property that it maps the constant field of $K$ to the constant field of $L$, and that this restriction to constant fields is a \emph{field} isomorphism (i.e., also additive). 

The method is now to prove compatibility of all constructions with field extensions, and then to embed residue fields as constant fields in large enough extensions of the original fields. 

So suppose that $\tilde{K}$ is an abelian extension of $K$, corresponding to a subgroup $H=\mathrm{Gal}(K^{\ab}/\tilde{K})$ of $G_K^{\ab}$, and let $\tilde{L}$ denote the abelian extension of $L$ corresponding to $\Phi(H)$. Here, we prove the analogue of Lemma 9 in \cite{Uchida}, namely, that the above constructed isomorphism $K^* \isomto L^*$ extends to an isomorphism $\tilde{K}^* \isomto \tilde{L}^*$. If $V_H$ denotes the transfer map (``Verlagerung'') from $G_K^{\ab}$ to $H$, we find a commutative diagram 
$$  \xymatrix{ K_v^* \ar@{->}[r]  \ar@{->}[d]  & \A_K^*  \ar@{->}[r]^{\rec_K} \ar@{->}[d]  & G_K^{\ab}\ar@{->}[d]^{V_H}   \\   \prod_{w\mid v} \tilde K_w^* \ar@{->}[r]  & \A_{\tilde K}^*  \ar@{->}[r]^{\rec_{\tilde K}} & H  },$$
in which the composed horizontal arrows are injective. 
From a similar diagram for the field $L$, and the fact that $\Phi \circ V_H = V_{\Phi(H)}$, we then find an induced isomorphism 
$$ \prod_{w\mid v} \tilde K_w^* \isomto \prod_{w'\mid \varphi(v)} \tilde L_{w'}^*, $$
and hence an isomorphism $ \A_{\tilde K}^* \isomto \A_{\tilde L}^*$ which is compatible with the reciprocity maps. Hence we find an isomorphism $\tilde K^* \isomto \tilde L^*$ which restricts to the constructed isomorphism $K^* \isomto L^*$. 

Now for any fixed place, there exists a finite extension $\tilde{K}$ of $K$ such that the residue field of $K_v$ is contained in the constant field of $\tilde{K}$. The above theory implies that the bijection of residue class fields $\bar K_v \rightarrow \bar L_{\varphi(v)}$ is the restriction of an isomorphism (of fields) between the constant fields of $\tilde K$ and $\tilde L$. This finishes the proof of additivity for residue fields. 

Finally, since the map $\Phi$ induces an isomorphism $\sigma$ of ground fields, we can twist one of the fields back by this Frobenius twist, to get a morphisms of function fields $K \cong L^{\sigma}$ fixing the ground field $k$, which then corresponds to an isomorphism of smooth curves. $\Box$

\section{Proof of Theorem \ref{main}} \label{??}

We write $K=k(X)$ and $L=k(Y)$, and we assume to have a degree-preserving isomorphism of dynamical systems $(S_X,I_X) \isomto (S_Y,I_Y)$. This means we have a degree preserving semigroup isomorphism $I_X \isomto I_Y$, hence a bijection of places 
$$ \varphi \colon \Places(K) \isomto \Places(L) \mbox{ with } \deg(\varphi(v))=\deg(v); $$
and an homeomorphism 
$ \Phi \colon S_X \isomto S_Y, $
which is equivariant with respect to the above map, i.e., 
$$ \Phi(v \ast x) = \varphi(v) \ast \Phi(x) \mbox{ for all } x \in S_X \mbox{ and } v \in \Places(K). $$ 
The first claim is: 

\begin{lemma} \label{lemgal} If there is a degree-preserving isomorphism $(\Phi, \varphi)$ of dynamical systems $(S_X,I_X) \isomto (S_Y,I_Y)$, then there is such an isomorphism for which 
the homeomorphism $\Phi$ induces a group isomorphism 
$ \Phi \colon G_{K}^{\ab} \isomto G_{L}^{\ab}. $ 
\end{lemma}

\noindent \emph{Proof.}\ 
The proof is inspired by the one in \cite{CM}, 7.3, but some differences arise from the geometric situation.  The  rough idea is as follows: a priori $\Phi$ is only a homeomorphism (between Cantor sets). However, using the equivariance of the isomorphism w.r.t.\ the action of ideals will allow us to prove that $\Phi$ is almost multiplicative ``on Frobenius elements for primes $\fp$'', which are dense in the $\fp$-unramified part of the \emph{geometric} Galois group $G^{\ab}_{\bar k K}$; one then extends to the entire Galois group. ``Almost multiplicative'' means that $\Phi \Phi(1)^{-1}$ is multiplicative on the geometric part, so making the resulting map genuinely multiplicative will involve changing it slightly. 

We take an integral ideal $\fm \in I_X$, and define $1_{\fm}$ to be the integral adele which is $1$ at the prime divisors of $\fm$ and zero elsewhere. Consider the subspace 
$H_{X,\fm}:=G_{K}^{\ab} \times_{\hat \cO_{K}^*} \{1_{\fm}\} \subseteq S_{X}.$
One immediately checks that 
$$ H_{X,\fm} \cong G_{K}^{\ab} / \rec_{X}(\prod_{\q \nmid \fm} \hat\cO_{\q}^*) \colon [(\gamma,1_{\fm})] \mapsto [\gamma] $$
is an isomorphism, and by class field theory, the right hand side is isomorphic to the Galois group of the maximal abelian extension of $K$ that is unramified \emph{outside} prime divisors of $\fm$. Consider its subgroup $G_{X,\fm}$, the Galois group of the maximal geometric abelian extension of $K$ unramified outside prime divisors of $\fm$. It corresponds to the subgroup $\bar H_{X,\fm} =  G_{\bar k K}^{\ab} \times_{\hat \cO_{K}^*} \{1_{\fm}\}$ of $H_{X,\fm}$. Now $G_{X,\fm}$ has a dense subgroup generated by $\rec_{X}(\fn)$ for $\fn$ running through the ideals $\fn$ that are coprime to $\fm$. Via the above map, this means that $\bar H_{X,\fm}$ is generated (as a group) by $\gamma_{\fn}:=[(\rec_{X}(\fn)^{-1},1_{\fm})]$ for $\fn$ running through the ideals coprime to $\fm$. Write  $\one_{\fm}=[(1,1_{\fm})]$, and $\Phi(\one_{\fm}) = [(x_{\fm},y_{\fm})]$. Since $\fm$ and $\fn$ are coprime, we have $[(\rec_{X}(\fn)^{-1},1_{\fm})] = [(\rec_{X}(\fn)^{-1},\fn 1_{\fm})]$, and hence we can write  $\gamma_{\fn} = \fn \ast \one_{\fm}$. In this way, we have written the generators in terms of the action. Now for two ideals $\fn_1$ and $\fn_2$ coprime to $\fm$: \begin{align*}
&\Phi(\one_{\fm}) \cdot  \Phi(\gamma_{\fn_1}   \cdot \gamma_{\fn_2}) \\ &= \Phi(\one_{\fm}) \cdot \left( \varphi(\fn_1) \varphi(\fn_2) \ast \Phi(\one_{\fm}) \right) \\ &= [(\rec_{Y}(\varphi(\fn_1) \varphi(\fn_2))^{-1} x_{\fm}^2, \varphi(\fn_1) \varphi(\fn_2) y_{\fm}^2)] \\ &= 
 [(\rec_{Y}(\varphi(\fn_1))^{-1} x_{\fm}, \varphi(\fn_1) y_{\fm})] \\  & \ \ \ \ \ \cdot  [(\rec_{Y}(\varphi(\fn_2))^{-1} x_{\fm}, \varphi(\fn_2) y_{\fm})] \\ &=  \left( \varphi(\fn_1) \ast \Phi(\one_{\fm}) \right) \cdot  \left( \varphi(\fn_2) \ast \Phi(\one_{\fm}) \right)  \\ &= \Phi(\gamma_{\fn_1}) \cdot \Phi(\gamma_{\fn_2}). 
\end{align*}
By density, we find that for all $\gamma_1, \gamma_2 \in \bar H_{X,\fm}$, we have 
$ \Phi(\one_{\fm}) \Phi(\gamma_{1} \gamma_{2})  = \Phi(\gamma_{1}) \Phi(\gamma_{2}). $ 

We compute the image of $\bar H_{X,\fm}$ under $\Phi$. Since $$H_{X,\fm} \subseteq \displaystyle \bigcap_{(\fm,\fn)=1} \fn \ast S_X,$$ the compatibility of actions implies that $$\Phi(H_{X,\fm}) \subseteq  \displaystyle \bigcap_{(\varphi(\fm),\fn')=1} \fn' \ast S_Y.$$ Choosing $\fn$ coprime to $\fm$, $\varphi(\fn)$ is coprime to $\varphi(\fm)$, and this inclusion shows that $y_{\fm}$ is zero on the support of $\varphi(\fn)$. Hence 
$ \Phi(\gamma_{\fn}) = [(\rec_{Y}(\varphi(\fn))^{-1} x_{\fm}, \varphi(\fn) y_{\fm})] = [(\rec_{Y}(\varphi(\fn))^{-1} x_{\fm}, y_{\fm})] \in \Phi(\one_{\fm}) G_{L}^{\ab}. $
By density, we conclude that $\Phi(\bar H_{K,\fm}) = \Phi(\one_{\fm}) G_{\bar k L}^{\ab}. $

Finally, we take direct limits under inclusion of ideals $\fm$. Since $\one_{\fm} \rightarrow 1$, we have $H_{X,\fm} \rightarrow G_{ K}^{\ab}$, and by continuity of $\Phi$ we find
$ \Phi \colon G_{K}^{\ab} \rightarrow \Phi(1) G_{L}^{\ab}$ has the property that for elements of $\gamma_1, \gamma_2 \in G_{\bar k K}^{\ab}$, it holds true that 
$ \Phi(1) \Phi(\gamma_{1} \gamma_{2})  = \Phi(\gamma_{1}) \Phi(\gamma_{2}). $ Now $\Phi(1)$ is not a zero divisor, since $$1 \in S_X \setminus \bigcup_{\fn \neq (1)} \fn \ast S_X$$ and equivariance implies that $$\Phi(1) \in S_Y \setminus \bigcup_{\fm \neq (1)} \fm \ast S_Y.$$ However, if $\Phi(1)=[(\gamma_0,\rho_0)]$ would be a zero divisor, some component $\rho_{0,\mathfrak{q}}$ would have to be zero, hence $\Phi(1) \in \mathfrak{q}\ast S_Y$.  
Thus, $$\tilde{\Phi}(\gamma):=\Phi(\gamma) \cdot \Phi(1)^{-1}$$ is a bijection $G_{K}^{\ab} \rightarrow G_{L}^{\ab}$ that restricts to a group isomorphism of geometric abelian Galois groups $G_{\bar k K}^{\ab} \rightarrow G_{\bar k L}^{\ab}$.  

Let $\bar S_X$ denote the ``geometric'' space $\bar S_X:=G_{\bar k K}^{\ab} \times_{\cO^*_X} \cO_X.$
There is a split exact sequence of topological groups 
$$ 0 \rightarrow G_{\bar k K}^{\ab} \rightarrow G_K^{\ab} \rightarrow G_k \cong \hat{\Z} \rightarrow 0,$$ by which we can fix an isomorphism $$G_{\bar k K}^{\ab} \times G_k \overset{i_K \times s_K}{\longrightarrow} G_K^{\ab}.$$ We claim that $S_X \isomto \bar S_X \times G_k$ is a homeomorphism by the map $[(\gamma,\rho)] \mapsto ([(\gamma_1,\rho)],\gamma_2)$ if we uniquely write $\gamma=\gamma_1 \gamma_2$ with $\gamma_1 \in i_K(G_{\bar k K}^{\ab})$ and $\gamma_2 \in s_K(G_k)$.  We replace the original homeomorphism $\Phi$ by the map $\Phi^\sim \colon S_X \isomto S_Y$ defined as 
$$S_X \cong \bar S_X \times G_k \overset{\frac{\Phi}{\Phi(1)} \times \mathrm{Id}}{\longrightarrow} \bar S_Y \times G_k \cong S_Y,$$ which changes nothing in the hypothesis; one may check by direct computation that $\Phi$ still intertwines the two actions of $I_X$ and $I_Y$ via the semigroup isomorphism $\varphi$. Hence we now assume without loss of generality that $\Phi$ restricts to a group isomorphism $G_K^{\ab} \rightarrow G_{L}^{\ab}$.

\begin{lemma} \label{lemlocal}
For any place $v$ of $K$, there is an isomorphism $\Phi_v$ of multiplicative groups of the corresponding local fields $ K_v^* \overset{\Phi_v}{\rightarrow} L_{\varphi(v)}^*$  with the property that $\ord_{\varphi(v)}(\Phi_v(x)) = \ord_v(x) \mbox{ for all }x \in K_v,$, and such that the corresponding diagram in Theorem \ref{thmuch} commutes. 
\end{lemma}

\noindent \emph{Proof.}\ 
We have a (non-canonically) split exact sequence $$0 \rightarrow \hat \cO_{v}^* \rightarrow K_v^* \rightarrow \langle \pi_v \rangle \rightarrow 0,$$ where $\pi_v$ is a uniformizer of $v$ and $\cO_v^*$ is the group of local integral units. By local class field theory, the latter is canonically identified with the inertia group of $v$ in $K^{\ab}$. 

We have seen in the previous proof that $\Phi$ respects ramification (and respects by construction constant field extensions), mapping $G_{X,\fm}$ to $G_{Y,\varphi(\fm)}$ (the Galois groups of the maximal extensions unramified outside $\fm$, respectively $\varphi(\fm)$). Basic Galois theory (as in Prop.\ 8.1 in \cite{CM}) then shows that  it also maps the group of the maximal abelian extension of $K$ unramified outside $\fm$ to that of $L$ ramified outside $\varphi(\fm)$, i.e., the inertia group of $v$ in $K^{\ab}$ is mapped to that of of $\varphi(v)$ in $L^{\ab}$. If we combine this with the map $\pi_v \mapsto \pi_{\varphi(v)}$, we find the desired map $\Phi_v \colon K_v^* \rightarrow L^*_{\varphi(v)}$ which automatically respects the order functions. The map $\Phi_v$ is also compatible with (local) reciprocity by construction, and the decomposition of the global reciprocity map as product of the local ones shows the product $\prod \Phi_v$ is compatible with the global reciprocity map, too. $\Box$

Theorem \ref{main} now follows by taking Lemma \ref{lemlocal} and Lemma \ref{lemgal} as input for Theorem \ref{thmuch}.

\section{Geometric version of Theorem \ref{cor}} \label{eff}

We recall some properties of $L$-series (cf.\ Chapter 9 in \cite{Rosen}). All $L$-series of $X$ are rational functions in $T=q^{-s}$. We have a (non-canonical) splitting of character groups $G_{k(X)}^{\ab \vee}  \cong G_{\bar k (X)}^{\ab \vee} \times \Z$. Denote a decomposition of a given character $\chi \in G^{\ab \vee}_{k(X)}$ according to this isomorphisms as $\chi = \chi_g \cdot \chi_c$, with $\chi_g \in G^{\ab \vee}_{\bar k (X)}$ and $\chi_c \in \Z$. Let $\mathrm{F}_k$ denote the Frobenius automorphism of the ground field $k$. If $\fp$ is a prime, then $\mathrm{Frob}_{\fp}$ acts like $F_{k}^{\deg \fp}$ in a constant field extension (\cite{Rosen}, 9.19), which proves the identity of Euler factors $$1-\chi(\mathrm{Frob}_{\fp}) T^{\deg \fp}=1-\chi_g(\mathrm{Frob}_{\fp}) \chi_c(\mathrm{F}^{\deg \fp}_k) T^{\deg \fp}.$$ Hence by expanding in an Euler product, we find 
$$ L_{k(X)}(\chi,T)   = L_{k(X)}(\chi_g, T \chi_c(\mathrm{F}_k)). $$
For the trivial character, this is the zeta function of $X$, which, when multiplied with $(1-T)(1-qT)$ is a polynomial of degree $2g$. For a character corresponding to a non-trivial constant extension $\F_{q^d}/\F_q$, it is equal to $\zeta_X(\chi(\mathrm{F}_k) T)$, which has poles, but not at $T=1$. If the character $\chi$ does not correspond to a constant extension of $X$, it is actually a polynomial in $T$ of degree $2g-2+\deg(\ff_\chi)$, where $g$ is the genus of $X$ and $\ff_\chi$ is the conductor of $\chi$ (\cite{Rosen}, 9.24A). 

\begin{theorem} \label{prop}
Let $X,Y$ denote two smooth projective curves over a finite field $k$. Then the following statements are equivalent: 
\begin{enumerate}
\item[\textup{(a)}] There is an isomorphism $\psi \colon G_{k(X)}^{\ab} \isomto G_{k(Y)}^{\ab}$ of Galois groups of the maximal abelian extensions of their respective function fields,  such that the corresponding $L$-series match: 
$ L(X,\chi) = L(Y,(\psi^{-1})^{\vee}\chi) $ for all characters $\chi \in G_{k(X)}^{{\ab} \vee}$. 
\item[\textup{(b)}] There is an isomorphism $\psi \colon G_{\bar k(X)}^{\ab} \isomto G_{\bar k(Y)}^{\ab}$ of Galois groups of the maximal geometric abelian extensions of their respective function fields,  such that the corresponding $L$-series match: 
$ L(X,\chi) = L(Y,(\psi^{-1})^{\vee}\chi) $ for all characters $\chi \in G_{\bar k(X)}^{{\ab} \vee}$. 
\end{enumerate}
\end{theorem}

\noindent \emph{Proof.}\ 
To prove that (a) implies (b), it suffices to show that a geometric character is mapped to a geometric character, which (by the  splitting) is equivalent to showing that a character corresponding to a constant extension is mapped to a similar character. For $\chi$ a character corresponding to a constant extension $\F_{q^d}/\F_q$, $L_{k(X)}(\chi,T) = \zeta_K(T \chi(\mathrm{F}_k))$ has a pole at $T$ equal to a primitive $d$-th root of unity $\chi_d(\mathrm{F}_k)^{-1}$. If it also equals an $L$-series $L_{k(Y)}(\chi',T)$ for some character $\chi'$, then $\chi'$ is forcedly the character of a constant extension (since there are poles) of degree $d$ (since there is a pole at a primitive $d$-th root of unity). This proves one direction.  In the other direction, we use the splitting to extend the isomorphism $\psi \colon G_{\bar k(X)}^{\ab} \isomto G_{\bar k(Y)}^{\ab}$ by the identity map on $\Z \rightarrow \Z$ to an isomorphism $\psi \colon G_{k(X)}^{\ab} \isomto G_{k(Y)}^{\ab}$. Then for any character $\chi \in G_{k(X)}^{\ab}$, we have \begin{align*} & L_{k(X)}(\chi,T) = L_{k(X)}(\chi_g, \chi_c(\mathrm{F}_k)T) \\ & = L_{k(Y)} ( \psi(\chi_g), \chi_c(\mathrm{F}_k)T) = L_{k(Y)} (\psi(\chi),T) \end{align*}
(where the middle equality arises from the assumption of matching of geometric $L$-series), as desired. $\Box$


\begin{thebibliography}{10}

\bibitem{BKT}
Fedor Bogomolov, Mikhail Korotiaev, and Yuri Tschinkel, \emph{A {T}orelli
  theorem for curves over finite fields}, Pure Appl. Math. Q. \textbf{6}
  (2010), no.~1, Special Issue: In honor of John Tate. Part 2, 245--294.

\bibitem{Magma}
Wieb Bosma, John Cannon, and Catherine Playoust, \emph{The {M}agma algebra
  system. {I}. {T}he user language}, J. Symbolic Comput. \textbf{24} (1997),
  no.~3-4, 235--265, Computational algebra and number theory (London, 1993).

\bibitem{BC}
Jean-Beno{\^{\i}}t Bost and Alain Connes, \emph{Hecke algebras, type {III}
  factors and phase transitions with spontaneous symmetry breaking in number
  theory}, Selecta Math. (N.S.) \textbf{1} (1995), no.~3, 411--457.

\bibitem{ConsMar}
Caterina Consani and Matilde Marcolli, \emph{Quantum statistical mechanics over
  function fields}, J. Number Theory \textbf{123} (2007), no.~2, 487--528.

\bibitem{CdJ}
Gunther Cornelissen and Jan~Willem de~Jong, \emph{The spectral length of a map
  between {R}iemannian manifolds}, preprint arxiv:1007.0907, to appear in J.\
  Noncommut.\ Geom.\, 2011.

\bibitem{CKZ}
Gunther Cornelissen, Aristides Kontogeorgis, and Lotte van~der Zalm,
  \emph{Arithmetic equivalence for function fields, the {G}oss zeta function
  and a generalisation}, J. Number Theory \textbf{130} (2010), no.~4,
  1000--1012.

\bibitem{CM}
Gunther Cornelissen and Matilde Marcolli, \emph{Quantum statistical mechanics,
  {$L$}-series, and anabelian geometry}, preprint arxiv:1009.0736, 2010.

\bibitem{DK}
Kenneth~R. Davidson and Elias~G. Katsoulis, \emph{Nonself-adjoint operator
  algebras for dynamical systems}, Operator structures and dynamical systems,
  Contemp. Math., vol. 503, Amer. Math. Soc., Providence, RI, 2009, pp.~39--51.

\bibitem{Gross}
Benedict~H. Gross, \emph{Trivial {$L$}-functions for the rational function
  field}, J. Number Theory (Carlitz memorial volume), in press, doi:
  10.1016/j.jnt.2012.02.011, 2012.

\bibitem{HP}
Eugene Ha and Fr{\'e}d{\'e}ric Paugam, \emph{Bost-{C}onnes-{M}arcolli systems
  for {S}himura varieties. {I}. {D}efinitions and formal analytic properties},
  IMRP Int. Math. Res. Pap. (2005), no.~5, 237--286.

\bibitem{Howe}
Everett~W. Howe, \emph{Constructing distinct curves with isomorphic
  {J}acobians}, J. Number Theory \textbf{56} (1996), no.~2, 381--390.

\bibitem{Jacob}
Beno{\^{\i}}t Jacob, \emph{Bost-{C}onnes type systems for function fields}, J.
  Noncommut. Geom. \textbf{1} (2007), no.~2, 141--211, post-publication
  corrected version: arxiv:math/0602554 (2012).

\bibitem{LLN}
Marcelo Laca, Nadia~S. Larsen, and Sergey Neshveyev, \emph{On {B}ost-{C}onnes
  types systems for number fields}, J. Number Theory \textbf{129} (2009),
  no.~2, 325--338.

\bibitem{NR}
Sergey Neshveyev and Simen Rustad, \emph{Bost-{C}onnes systems associated with
  function fields}, preprint, arxiv:1112.5826, 2011.

\bibitem{Rosen}
Michael Rosen, \emph{Number theory in function fields}, Graduate Texts in
  Mathematics, vol. 210, Springer-Verlag, New York, 2002.

\bibitem{Tate}
John Tate, \emph{Endomorphisms of abelian varieties over finite fields},
  Invent. Math. \textbf{2} (1966), 134--144.

\bibitem{Turner}
Stuart Turner, \emph{Adele rings of global field of positive characteristic},
  Bol. Soc. Brasil. Mat. \textbf{9} (1978), no.~1, 89--95.

\bibitem{Uchida}
K{\^o}ji Uchida, \emph{Isomorphisms of {G}alois groups of algebraic function
  fields}, Ann. Math. (2) \textbf{106} (1977), no.~3, 589--598.

\bibitem{Zilber}
Boris Zilber, \emph{A curve and its abstract {J}acobian}, preprint, available
  from people.maths.ox.ac.uk/zilber, 2012.

\end{thebibliography}
\end{document}